\input jytex.tex   
\typesize=10pt
\magnification=1200
\baselineskip17truept
\hsize=6truein\vsize=8.5truein
\sectionnumstyle{blank}
\chapternumstyle{blank}
\chapternum=1
\sectionnum=1
\pagenum=0

\def\begintitle{\pagenumstyle{blank}\parindent=0pt\begin{narrow}[0.4in]}
\def\endtitle{\end{narrow}\newpage\pagenumstyle{arabic}}


\def\beginexercise{\vskip 20truept\parindent=0pt\begin{narrow}[10
truept]}
\def\endexercise{\vskip 10truept\end{narrow}}


\def\eql#1{\eqno\eqnlabel{#1}}
\def\ref{\reference}
\def\peq{\puteqn}
\def\pref{\putref}

\def\mgn{\marginnote}
\def\bex{\begin{exercise}}
\def\eex{\end{exercise}}


\font\open=msbm10 
\def\mbox#1{{\leavevmode\hbox{#1}}}

\def\hspace#1{{\phantom{\mbox#1}}}
\def\oZ{\mbox{\open\char90}}

\def\be{\beta}
\def\ga{\gamma}
\def\de{\delta}
\def\Ga{\Gamma}

\def\ka{\kappa}
\def\la{\lambda}
\def\La{\Lambda}
\def\om{\omega}

\def\Si{\Sigma}
\def\th{\theta}

\def\ze{\zeta}

\def\De{\Delta}

\def\caE{{\cal E}}
\def\caL{{\cal L}}

\def\Real{{\rm Re\,}}

\def\zf{$\zeta$--function}
\def\zfs{$\zeta$--functions}
\def\hk{heat-kernel}
\def\ol{\overline}


\def\frac#1/#2{\leavevmode\kern.1em
\raise.5ex\hbox{\the\scriptfont0 #1}\kern-.1em/\kern-.15em
\lower.25ex\hbox{\the\scriptfont0 #2}}
\def\sfrac#1/#2{\leavevmode\kern.1em
\raise.5ex\hbox{\the\scriptscriptfont0 #1}\kern-.1em/\kern-.15em
\lower.25ex\hbox{\the\scriptscriptfont0 #2}}

\def\gtorder{\mathrel{\raise.3ex\hbox{$>$}\mkern-14mu
             \lower0.6ex\hbox{$\sim$}}}
\def\ltorder{\mathrel{\raise.3ex\hbox{$<$}\mkern-14mu
             \lower0.6ex\hbox{$\sim$}}}

\def\semidirprod{\rlap{\ss C}\raise1pt\hbox{$\mkern.75mu\times$}}
\def\for{\lower6pt\hbox{$\Big|$}}
\def\fish{\kern-.25em{\phantom{abcde}\over \phantom{abcde}}\kern-.25em}


\def\boxit#1{\vbox{\hrule\hbox{\vrule\kern3pt
        \vbox{\kern3pt#1\kern3pt}\kern3pt\vrule}\hrule}}
\def\dalemb#1#2{{\vbox{\hrule height .#2pt
        \hbox{\vrule width.#2pt height#1pt \kern#1pt
                \vrule width.#2pt}
        \hrule height.#2pt}}}
\def\square{\mathord{\dalemb{5.9}{6}\hbox{\hskip1pt}}}

\def\frac#1#2{{{#1}\over{#2}}}


\def\etc{{\it etc. }}

\def\eg{{\it e.g. }}
\def\ie{{\it i.e. }}
\def\cf{{\it cf }}
\def\pa{\partial}


\def\Tr{{\rm Tr\,}}

\def\3j#1#2#3#4#5#6{\left\lgroup\matrix{#1&#2&#3\cr#4&#5&#6\cr}
\right\rgroup}

\def\man{{\cal M}}
\def\caI{{\cal I}}

\def\can{{\cal N}}

\def\m?{\mgn{?}}

\def\pa{\partial}

\def\beq{\begin{eqnarray}}
\def\eeq{\end{eqnarray}}


\def\aop#1#2#3{{\it Ann. Phys.} {\bf {#1}} ({#2}) #3}

\def\cmp#1#2#3{{\it Comm. Math. Phys.} {\bf {#1}} ({#2}) #3}
\def\cqg#1#2#3{{\it Class. Quant. Grav.} {\bf {#1}} ({#2}) #3}

\def\jgp#1#2#3{{\it J. Geom. and Phys.} {\bf {#1}} ({#2}) #3}
\def\jmp#1#2#3{{\it J. Math. Phys.} {\bf {#1}} ({#2}) #3}
\def\jpa#1#2#3{{\it J. Phys.} {\bf A{#1}} ({#2}) #3}

\def\np#1#2#3{{\it Nucl. Phys.} {\bf B{#1}} ({#2}) #3}
\def\pl#1#2#3{{\it Phys. Lett.} {\bf {#1}} ({#2}) #3}

\def\pr#1#2#3{{\it Phys. Rev.} {\bf {#1}} ({#2}) #3}

\def\prD#1#2#3{{\it Phys. Rev.} {\bf D{#1}} ({#2}) #3}
\def\prE#1#2#3{{\it Phys. Rev.} {\bf E{#1}} ({#2}) #3}

\def\cras#1#2#3{{\it Comptes Rend. Acad. Sci. (Paris)} {\bf{#1}} (#2) #3}

\def\mpcps#1#2#3{{\it Math. Proc. Camb. Phil. Soc.} {\bf{#1}} ({#2}) #3}

\def\am#1#2#3{{\it Acta Mathematica} {\bf {#1}} ({#2}) #3}
\def\aim#1#2#3{{\it Adv. in Math.} {\bf {#1}} ({#2}) #3}
\def\ajm#1#2#3{{\it Am. J. Math.} {\bf {#1}} ({#2}) #3}
\def\amm#1#2#3{{\it Am. Math. Mon.} {\bf {#1}} ({#2}) #3}

\def\aom#1#2#3{{\it Ann. of Math.} {\bf {#1}} ({#2}) #3}

\def\cpde#1#2#3{{\it Comm. Partial Diff. Equns.} {\bf {#1}} (19{#2}) #3}

\def\invm#1#2#3{{\it Invent. Math.} {\bf {#1}} ({#2}) #3}
\def\ijpam#1#2#3{{\it Ind. J. Pure and Appl. Math.} {\bf {#1}} (19{#2}) #3}
\def\jdg#1#2#3{{\it J. Diff. Geom.} {\bf {#1}} ({#2}) #3}
\def\jfa#1#2#3{{\it J. Func. Anal.} {\bf {#1}} ({#2}) #3}

\def\jmpa#1#2#3{{\it J. Math. Pures. Appl.} {\bf {#1}} ({#2}) #3}

\def\ojm#1#2#3{{\it Osaka J.Math.} {\bf {#1}} ({#2}) #3}

\def\pja#1#2#3{{\it Proc. Jap. Acad.} {\bf {A#1}} ({#2}) #3}

\def\tams#1#2#3{{\it Trans. Am. Math. Soc.} {\bf {#1}} ({#2}) #3}

\begin{title}
\vglue 1truein
\vskip15truept
\centertext {\Bigfonts \bf The hybrid spectral problem and} \vskip5truept
\centertext {\Bigfonts \bf  Robin boundary conditions}

\vskip 20truept \centertext{J.S.Dowker\footnote{dowker@a35.ph.man.ac.uk}}
\vskip 7truept \centertext{\it School of Physics and Astronomy, }
\centertext{\it Theoretical Physics Group, } \centertext{The University of
Manchester,} \centertext{ Manchester, England} \vskip 20truept \centertext
{Abstract} \vskip10truept
\begin{narrow}
The hybrid spectral problem where the field satisfies Dirichlet conditions
(D) on part of the boundary of the relevant domain and Neumann (N) on the
remainder is discussed in simple terms. A conjecture for the $C_1$
coefficient is presented and the conformal determinant on a 2-disc, where
the $D$ and $N$ regions are semi-circles, is derived. Comments on higher
coefficients are made.

A separable second order hemisphere hybrid problem is introduced that
involves Robin boundary conditions and leads to logarithmic terms in the
heat--kernel expansion which are evaluated explicitly.
\end{narrow}
\vskip 5truept
\vskip 60truept
\vfil
\end{title}
\pagenum=0
\section{\bf1. Introduction}
The explicit construction of the general form of the \hk\ expansion
coefficients has reached the stage when further progress is impeded mainly
by ungainliness. Unless there is some compelling reason for finding a
specific higher coefficient, its exhibition is not particularly enlightening
and is not really worth the, often considerable, effort. Other, more
productive, avenues consist of generalising the differential operator, the
manifold or the boundary conditions. In the latter context a simply stated
extension is the class of problems where the field satisfies Dirichlet
conditions (D), say, on part of the boundary and Neumann (N) on the
remainder. These boundary conditions are sometimes termed `mixed ' in the
classical literature (\eg Sneddon [\pref{Sneddon}]) or sometimes `hybrid'
(\eg Treves [\pref{treves2}] chap.37). A brief history of the corresponding
potential theory, sometimes referred to as the Zaremba problem, is contained
in Azzam and Kreysig, [\pref{AandK}]. It is also interesting to note that
these conditions have occurred in string theory, [\pref{Siegel}], and
have recently been considered in connection with isospectrality, Jakobson
{\it et al}, [\pref{JLNP}].

To set some notation, the conventional short--time expansion of the
integrated heat--kernel of a smooth boundary value problem is,
  $$
  K(t)\equiv\Tr e^{-Pt}\sim{1\over(4\pi t)^{d/2}}
  \sum_{n=0,1/2,1\ldots}^\infty C_n t^n\,,
  \eql{hkasymp}
  $$
where $P$ is a smooth (singularity free) elliptic differential operator and
where, initially, the manifold, its boundary and any boundary conditions are
all smooth. Typically $P$ is the Laplacian, plus possibly a smooth
potential, and the coefficients are locally computed as integrals over the
manifold, or its boundary, of local geometrical invariants constructed from
the curvature, for example. See [\pref{Gilkey1}] for an extensive treatment.

Any relaxation of smoothness can result in modifications to this expansion.
For example, a singular potential can lead to `anomalous' powers of $t$, \eg
[\pref{FMPS}].

The hybrid N/D case has non--smooth boundary conditions and so can be
classed as a singular boundary value problem, see \eg [\pref{avramidi,
avramidi2}]. Even though it may be geometrically smooth, the codimension-2
region, $\Si$, where the D and N conditions meet, can be regarded as a
conical singularity. The general existence of logarithmic terms in the
expansion of $K(t)$ for such singularities, and other situations, has been
analysed for a long time (\eg Cheeger [\pref{cheeg1}], Br\"uning and Seeley
[\pref{BandS}], Grubb and Seeley [\pref{GrandS}], Grubb [\pref{Grubb,
Grubb4}], Gilkey and Grubb [\pref{GandG}], Seeley [\pref{Seeley2}]).

Seeley, [\pref{Seeley}], based on [\pref{BandS}], has proved the existence
of an asymptotic expansion that allows for logarithmic terms but then shows
for a particular N/D situation that such terms do not appear. A classic
case is the simple $\pi$--wedge where the explicit calculation shows that
logarithmic terms are absent.

It is anticipated that the \hk\ coefficients will receive contributions from
$\Si$. This has been confirmed by Avramidi [\pref{avramidi,avramidi2}] and
work by  van den Berg and Gilkey, [\pref{VdBandG}], on heat {\it content} is
also pertinent.

In this paper I wish to discuss some aspects of the hybrid question that are
mainly example driven and with a minimum of algebra. It is hoped that these
considerations will prove useful in more general field and string theoretic
areas where \hk\ coefficients play important roles in divergence and scaling
questions. My approach is mostly global, as opposed to the local treatment
by Avramidi, [\pref{avramidi, avramidi2}].

I begin with the Laplacian eigenproblem on the interval with D and N
conditions on the ends. This is then embedded in higher dimensions and used
to determine the $C_1$ hybrid coefficient for the 2--wedge from the purely
D, or purely N, expressions, which are very old. This result is then used to
write down the general $C_1$ for a $d$--manifold with piecewise smooth
boundary and a conjecture is made for the case where N is replaced by Robin
conditions, denoted (R). A crude check in the N--D case results from
applying the technique to the 2--lune.

Although  my main attention is directed at $C_1$, some very limited
information on the higher coefficients, $C_{3/2}$ and $C_2$, is extracted in
section 5 from the hybrid half--disc.

As an example of the use of the conjectured form of the hybrid $C_1$, I
evaluate the Laplacian functional determinant on the N--D disc in section 6.

In section 7, I set up a separable Robin hybrid problem for the Laplacian on
the 2--sphere and show in later sections that logarithmic terms appear in
the expansion of the heat--kernel. Perturbation theory is used to bolster
confidence in the existence of the model.

\section{\bf2. Basic idea}

For the Laplacian, simple calculation, or the drawing of a few modes, shows
that on the interval of length $L$ with Dirichlet (D) and Neumann (N)
conditions, the spectral data of the various eigenproblems are related by,
 $$\eqalign{ (D,N)_L\cup(D,D)_L&=(D,D)_{2L}\cr
(D,N)_L\cup(N,N)_L&=(N,N)_{2L}\cr}
\eql{relns}
 $$
 $$
(D,D)_L\cup(N,N)_L=P_{2L}\,,
 \eql{relns2}
 $$
where the notation $(D,N)$ signifies a problem with $D$ conditions at one
end and $N$ at the other and $P$ stands for periodic conditions. Averaging
(\peq{relns}) gives, using (\peq{relns2})
  $$
(D,N)\cup{1\over2}P_{2L}={1\over2}P_{4L}\,.
 \eql{relns3}
 $$
These relations can be extended algebraically to any spectral quantity, such
as the heat--kernel.

The `subtraction' implied by (\peq{relns}) and (\peq{relns3}), in order to
extract the $(D,N)$ part, amounts to a cull of the even modes on the doubled
interval, as is well known (\cf Rayleigh [\pref{Rayleigh}], vol I, p.247).

The relations can be applied to the arc of a circle, which might form part
of an SO(2) foliation of a two--dimensional region (or the projection of a
higher dimensional region onto two dimensions). A wedge is a good example
which I will now look at. Separability of the modes implies that the
relations (\peq{relns}) apply equally well to the wedge, where the notation
signifies that either $D$ or $N$ holds on the straight sides, (say $\th=0$
and $\th=\be$). Equation (\peq{relns}) can then immediately be applied to
the heat--kernel, and its small--time expansion, to determine the form of
the heat--kernel coefficients in the $(D,N)$ combination. I will show how
this works out for the $C_1$ coefficient.

The $(D,D)$ and $(N,N)$ wedge coefficients are well known and have been
derived in several ways. They are,
  $$ C^{\rm wedge}_1(D,D)=C^{\rm wedge}_1(N,N)={\pi^2-\be^2\over6\be}\,,
\eql{cee1s}
  $$
Hence from (\peq{relns})
  $$ C^{\rm wedge}_1(D,N)=-{\pi^2+2\be^2\over12\be}.
\eql{cee1DN}
  $$
This last result has been derived by Watson in a rather complicated way
using the modes directly, [\pref{Watson}].

Incidentally the conjecture by Gottlieb (equn.(3.5) in [\pref{Gottlieb}]),
that the $(N,D)$ case differs from the $(D,D)$ one only by a sign, is
incorrect, although it is true in the special case of a right-angled wedge,
as is easily checked by looking at rectilinear flat domains. The statement
is carried through into ref [\pref{DandG}].

Sommerfeld, [\pref{Sommerfeld}] vol.2 p.827, also mentions the `mixed' wedge
and indicates how to treat it using images if $\be=\pi n/m$.

It is useful to note that, as pointed out by Cheeger, [\pref{cheeg1}] p.605,
the expressions (\peq{cee1s}) and (\peq{cee1DN}), are not locally computable
geometric invariants, as evidenced by the $1/\be$ dependence. In the
derivation of the wedge coefficients by Cheeger (see Bordag {\it et al},
[\pref{BKD}], Cognola and Zerbini, [\pref{CandZ}]) the $1/\be$ terms arise
from the $\be$--interval \zf\ evaluated at the argument $-1/2$. As noted in
[\pref{BKD}] this is the Casimir energy on the interval, clarifying  the
nonlocal character of this term. By contrast, the term proportional to $\be$
is locally computable.

It is an important technical point that, as discussed by Avramidi,
[\pref{avramidi,avramidi2}], and Seeley, [\pref{Seeley}], it is necessary to
specify boundary conditions at the singular region $\Si$ to give a well
defined problem. The derivations of the expressions (\peq{cee1s}) assume
Dirichlet at the wedge apex. This amounts to taking the Friedrichs extension
by default and extends to the hybrid wedge (\peq{cee1DN}) by (\peq{relns}).
I will continue with this simplifying choice for the rest of this paper.

\section{\bf3. The general case}
Consider in general dimension a manifold whose boundary is piecewise smooth
consisting of domains, $\pa\man_i$, which intersect in codimension-2
manifolds, ${\cal I}_{ij}$. On each of the pieces, $\pa\man_i$, either $D$
or $N$ is imposed. Kac's principle of not feeling the boundary (a cheap way
of avoiding estimates), [\pref{Kac}], implies that, $C_1$ will take
contributions from the manifolds of codimension zero, one and two
independently.

In general dimension, for all $D$ or all $N$, the smeared coefficients are
known,
  $$
 C_1(D)=\bigg({1\over6}-\xi\bigg)\int_{\man}\!\! Rf dV
 +\int_{\pa\man}\!\!\bigg({1\over3}\ka -{1\over2}\,n\cdot\pa\bigg)f\,dA
 +{1\over6}\int_{\cal I}\!\!{\pi^2-\be^2\over\be}\,f\,dL
 \eql{cee1d}
 $$
  $$
  C_1(N)=\bigg({1\over6}-\xi\bigg)\int_{\man}\!\! Rf dV
  +\int_{\pa\man}\!\!\bigg({1\over3}\ka
 + {1\over2}\,n\cdot\pa\bigg)f\,dA+
 {1\over6}\int_{\cal I}\!\!{\pi^2-\be^2\over\be}\,f\,dL
 \,,
 \eql{cee1n}
   $$
where $\pa\man$ is the union of the boundary pieces and ${\cal I}$ that of
the intersections,
  $$
  \pa\man=\bigcup_i\, \pa\man_i,\quad {\cal I}=\bigcup_{i<j}\,{\cal I}_{ij}\,.
  $$

The `smeared' coefficients result from the trace $\Tr fe^{-Pt}$, where $f$
is a spatially local operator, and are handy when discussing conformal
transformations and for extracting the integrands. I will not use this
freedom in any overt calculational way (see Kirsten, [\pref{Kirstenb}]).

For a mixture of $D$ and $N$, the volume contribution clearly remains
unchanged while the surface contribution divides simply into a sum
separately over those regions $\pa\man(D)$ and $\pa\man(N)$ subject to $D$
and $N$ respectively. The codimension--2 intersections ${\cal I}_{ij}$
divide into the three (wedge) types ${\cal I}(D,D)$, ${\cal I}(N,N)$ and
${\cal I}(N,D)$ so the corresponding $C_1$ is, using (\peq{cee1DN}),
 $$ \eqalign{C_1(D,N)&=\bigg({1\over6}-\xi\bigg)\int_\man Rf
dV+\int_{\pa\man(D)}\bigg({1\over3}\,\ka
-{1\over2}\,n\cdot\pa\bigg)f\,dA\cr&+\int_{\pa\man(N)}\bigg({1\over3}\ka +
{1\over2}\,n\cdot\pa\bigg)f\,dA+{1\over6}\int_{{\cal I}(D,D)\cup{\cal
I}(N,N)} {\pi^2-\be^2\over\be}\,f\,dL\cr &-{1\over12}\int_{{\cal
I}(D,N)}{\pi^2+2\be^2\over\be}f\,dL\,. \cr}\eql{ceedn1}
 $$

In accordance with a previous remark, the wedge-like codimension--2
contributions in the above expressions are not locally computable.

If the boundary is smooth, then all the dihedral angles $\be$ equal $\pi$
and the codimension--2 part of (\peq{ceedn1}) (the last two integrals)
reduces to

  $$ -{\pi\over4}\int_{{\cal I}(D,N)}f\,dL\,.
\eql{cee12}
  $$
To repeat, even though the boundary is smooth, the region ${\cal
I}(D,N)\equiv \Si$ is a singular region.

For example, for the 3-ball with $D$ on the northern hemisphere and
$N\,(S=0)$ on the southern,
  $$ C_1(D,N)={8\pi\over3}-{\pi^2\over2}\,,
  $$
for the smearing function, $f=1$.

A local derivation of (\peq{cee12}), justifying the application of Kac's
principle, has been given by Avramidi [\pref{avramidi,avramidi2}]. It has
also been obtained by van den Berg (unpublished).

For later reference, I would like to extend the Neumann conditions to Robin,
(R), ones,
  $$
  \big(n\cdot\pa-S\big)\Phi\bigg|_{\pa\man}=0\,,
  \eql{robin}
  $$
where $n$ is the inward normal and $S$ can depend on position. For example,
in place of (\peq{cee1n}) one might expect,
  $$
  C_1(R)=\bigg({1\over6}-\xi\bigg)\int_{\man}\!\! Rf dV
  +\int_{\pa\man}\!\!\bigg({1\over3}\ka
 -2S+ {1\over2}n\cdot\pa\bigg)f\,dA+
 {1\over6}\int_{\cal I}\!\!{\pi^2-\be^2\over\be}\,f\,dL
 \,.
 \eql{cee1r}
   $$

The first two terms are the standard ones, \eg [\pref{BandG2,Gilkey1}], for
a smooth boundary. The last, codimension--2, term has actually not been
derived directly for the Robin wedge but it holds when $S=0$ and dimensions
show that the Robin function cannot enter algebraically into this term.

On the same basis my conjecture for the corresponding $C_1(D,R)$ is,
 $$ \eqalign{C_1(D,R)&=\bigg({1\over6}-\xi\bigg)\int_\man Rf
dV+\int_{\pa\man(D)}\bigg({1\over3}\,\ka
-{1\over2}n\cdot\pa\bigg)f\,dA\cr&+\int_{\pa\man(R)}\bigg({1\over3}\ka -2S+
{1\over2}n\cdot\pa\bigg)f\,dA+{1\over6}\int_{{\cal I}(D,D)\cup{\cal I}(R,R)}
{\pi^2-\be^2\over\be}\,f\,dL\cr &-{1\over12}\int_{{\cal
I}(D,R)}{\pi^2+2\be^2\over\be}f\,dL\,. \cr} \eql{cee1}
 $$

Proceeding on the basis that (\peq{cee1}) is correct, an expression for
$C_1(N,R)$ can be obtained by making $\pa\man(D)$ empty and dividing
$\pa\man(R)$ into $\pa\man(N)\cup\pa\man(R)$. the conjecture is then,
  $$ \eqalign{
  C_1(N,R)&=\bigg({1\over6}-\xi\bigg)\int_\man Rf
 dV+\int_{\pa\man}\bigg({1\over3}\ka + {1\over2}n\cdot\pa\bigg)f\,dA\cr
 &-2\int_{\pa\man(R)}S f\,dA+{1\over6}\int_{{\cal I}}
 {\pi^2-\be^2\over\be}\,f\,dL\,.\cr
 }
 \eql{cee1nr}
 $$

It is clear that the above expressions have a specific validity, even
without the codimension--2 parts. Thus, although the Robin form
(\peq{cee1r}) trivially reduces to the Neumann one, (\peq{cee1n}), when
$S=0$, it is not possible to obtain the Dirichlet form, (\peq{cee1d}),
directly by setting $S=\infty$. The exhibited forms, which, to repeat, are
coefficients in a `small time' asymptotic expansion, are really valid in the
limit of small $S^2t$, as discussed by Fulling, [\pref{fulling}]. I return
more specifically to Robin conditions from section 7 onwards

\section{\bf4. The lune}

The expression for $C_1(D,N)$ can be checked in a curved space case by
considering a lune segment of a sphere.

Relations (\peq{relns}), (\peq{relns2}) can be applied to the 2-lune where
the intervals are the sections of the lines of latitude cut out by the the
two longitudes, $\phi=0$, $\phi=\be$. In this case the extrinsic curvatures
vanish (the boundaries are geodesically embedded) but there is a volume
(area) term independent of the boundary conditions.

The \zfs\ are now somewhat more explicit [\pref{D2,ChandD}]. It is possible
to work with general angle $\be$. If one chooses $\be=\pi/q$, $q\in\oZ$, the
\zfs\ have been derived in [\pref{ChandD}] and used in [\pref{AandD}].

Denoting the lune by $\caL(\be)$ one has,
$$\eqalign{
(D,N)_{\caL(\be)}\cup(D,D)_{\caL(\be)}&=(D,D)_{\caL(2\be)}\cr
(D,N)_{\caL(\be)}\cup(N,N)_{\caL(\be)}&=(N,N)_{\caL(2\be)}\,,\cr}
$$
so that the corresponding \zfs\ combine algebraically,
$$
\ze_\be^{ND}(s)=\ze_{2\be}^{DD}(s)-\ze_\be^{DD}(s)
=\ze_{2\be}^{NN}(s)-\ze_\be^{NN}(s)\,.
\eql{ndzeta}$$

The $DD$ and $NN$
\zfs\ have been derived in [\pref{ChandD}] as Barnes \zfs\ for conformal
coupling in three dimensions (leading to simple eigenvalues) and yield the
specific value, for example,
$$
\ze_\be^{DD}(0)={1\over12}\bigg({\pi\over\be}-{\be\over2\pi}\bigg),
$$
which can be used to confirm the expression (\peq{cee1DN}) using the relation
between $C_1$ and $\ze(0)$. (In this case there are no zero modes.)

The volume contribution, $\be/6$, to $C_1$ is standard and is the same for all
boundary conditions. Hence the contribution of each $(N,D)$ corner
(of which there are two) is
 $$
{1\over2}\bigg[-{4\pi\over24}\bigg({\pi\over\be}+{\be\over \pi}\bigg)
-{\be\over6}\bigg]=
-{\pi^2+2\be^2\over12\be}
 $$
as required for the check.
\section{\bf5. The disc and semi-circle. Higher coefficients.}
The fact that the extrinsic curvatures are zero means that the lune is not
excessively helpful in deriving the form of the higher coefficients in the
$(D,N)$ case. Some further information can be obtained by looking at the
half-disc  with semi-circular boundary having different conditions on the
diameter and circumference.

A straightforward application of, say the Stewartson and Waechter Laplace
transform technique combined with an image method soon yields the results
for the short time expansions
 $$
K_{DD}(t)\sim{1\over8t}-{2+\pi\over8\sqrt{\pi t}}+{5\over24}
+{\sqrt t(\pi+16)\over256
\sqrt\pi}+\bigg({1\over315}+{1\over32}\bigg)\,t+\ldots
\eql{ssDD}
 $$
 $$
K_{ND}(t)\sim{1\over8t}+{2-\pi\over8\sqrt{\pi t}}-{1\over24}
+{\sqrt t(\pi-16)\over256
\sqrt\pi}+\bigg({1\over315}-{1\over32}\bigg)\,t+\ldots
\eql{ssND}
 $$
 $$
K_{NN}(t)\sim{1\over8t}+{2+\pi\over8\sqrt{\pi t}}+{5\over24}
+{\sqrt t(5\pi+48)\over256
\sqrt\pi}+\bigg({1\over45}+{3\over32}\bigg)\,t+\ldots
\eql{ssNN}
 $$
 $$
K_{DN}(t)\sim{1\over8t}-{2-\pi\over8\sqrt{\pi t}}-{1\over24}
+{\sqrt t(5\pi-48)\over256
\sqrt\pi}+\bigg({1\over45}-{3\over32}\bigg)\,t+\ldots
 \eql{ssDN}$$
where $DN$ means $D$ on the diameter and $N$ on the circumference, \etc

The constant terms check with (\peq{cee1s}) and (\peq{cee1DN}) for
$\be=\pi/2$. Also (\peq{relns2}), applied to the diameter as a wedge of
angle $\pi$, yields the $D$ and $N$ (\eg [\pref{Moss}]), {\it full} disc
expansions.

The extrinsic curvature vanishes on the diameter and equals one on the
circumference part of the boundary so some information on the $C_{3/2}$ and
$C_2$ coefficients can be inferred. Formulae in the non-mixed types $(D,D)$
and $(N,N)$ have been given in [\pref{AandD,DandA}] which agree with the
relevant parts of the above expressions. Indeed I used the hemi-disc in
deriving these results.

 Also in [\pref{AandD}] will be found an expression for $C_2$ in the case
the boundary parts $\pa\man_i$ are subject to Robin conditions with
different boundary functions, $S_i$ although all dihedral angles are
restricted to $\pi/2$.

In the case of $C_2$, the $1/315$ is the contribution of the curved $D$
semicircle while the $\pm1/32$ is the effect of the (two) corners and
likewise regarding the $1/45\pm3/32$ combination. The $C_{3/2}$ coefficient
exhibits a similar structure. Experience with the flat wedge shows that it
is unwise to draw too many conclusions when the angle is $\pi/2$. What we
can say, however, is that, using the $3/2$ coefficient as an exemplar, one
term will have the general form
 $$\eqalign{
-{\sqrt\pi\over24}\bigg[\int_{\caI(D,D)}\la_{DD}&(\be)(\ka_1+\ka_2)\,dL+
\int_{\caI(N,N)}\la_{NN}(\be)(\ka_1+\ka_2)\,dL +\cr
&\int_{\caI(N,D)}\big(\la_{ND}(\be)\ka_D+\la_{DN}(\be)\ka_N\big)\,dL\bigg]\cr}
 $$
where $\la_{ND}(\pi/2)=-\la_{DD}(\pi/2)=3$ and $\la_{DN}(\pi/2)=
-\la_{NN}(\pi/2)=-9$.
This change of sign is a simple consequence of images, or of (\peq{relns})
since the
$DD$ and $NN$ quantities vanish when $\be=\pi$.
\section{\bf6. The disc determinant.}
A direct attack via modes, of what is, after rectilinear domains, the
simplest two-dimensional situation, \ie a disc subject to $N$ on one half of
the circumference, and $D$ on the rest, would seem to be difficult in so far
as the construction of the \zf\ or heat-kernel is concerned. However, the
functional determinant, defined conventionally  as $\exp\big(-\ze'(0)\big)$,
appears to be accessible by conformal transformation from that on an
$ND$-lune of angle $\pi$, \ie a hemisphere with $N$ on one half of the rim
(the equator) and $D$ on the rest, which is an easy quantity to find in
terms of Barnes \zf\ from (\peq{ndzeta}).

Instead of the determinant I use the effective action, $W$, defined by $W=
-\ze'(0)/2$. Integrating the conformal anomaly leads to the relation
  $$
   W[e^{-2\om} g]=W[g]+W[e^{-2\om} g,g]
  $$
where $W[e^{-2\om} g,g]$ is the cocyle function.

For the above programme to work, one would need the conjectured form of
$C_1$, (\peq{cee1}), to be valid in order to construct the required cocycle
function in two dimensions. Applying the standard techniques this is (\cf
[\pref{D}]), for a {\it smooth} boundary,
 $$\eqalign{
W[e^{-2\om} g,g]={1\over24\pi}\int_\man\om\big(R+\square\om\big)\,dV
+{1\over12\pi} \int_{\pa\man}\om\big(\ka+{1\over2}(n.\pa)\om\big)\,dA+\cr
{1\over8\pi}\bigg(\int_{\pa\man(N)}-\int_{\pa\man(D)}\bigg)(n.\pa)\,\om\,dA
-{1\over16}\sum_k\om_k\cr} \eql{cocycle}
 $$
where $k$ labels the points where $D$ and $N$ meet and $\om_k$ are the
values of $\om$ at these points. If $\pa\man(D)$ is empty there is a volume
term coming from the pure $N$ zero mode.

To go from the hemisphere to the disc I employ the equatorial stereographic
projection as in [\pref{Weis,BandG,AandD,DandA,D}] noting that there is no
codimension-2 contribution because the conformal factor is unity on the
boundary, implying $\om_k=0$.

Then (\peq{cocycle}) can be written
 $$
W_{ND}[\bar g,g]={1\over2}\bigg(W_{D}[\bar g,g]
+\overline W_{N}[\bar g,g]\bigg),
\eql{cocycle2}
 $$
where $\overline W_{N}$ means the usual Neumann expression, omitting the
zero mode piece, and I can use the known values, ($\bar g$ = disc and $g$ =
hemisphere),
 $$\eqalign{
W_D[\bar g,g]&={1\over6}\log2-{1\over3}\cr
\overline W_N[\bar g,g]&={2\over3}\log2+{1\over6}.\cr}
\eql{cocycle3}
 $$

The \zf\ on the $ND$-hemisphere follows from (\peq{ndzeta}) with $\be=\pi$.
The \zf, $\ze_\pi^{DD}(s)$ is the usual hemisphere \zf\ and the determinant has
been considered a number of times.
$\ze_{2\pi}^{DD}(s)$, corresponds to Sommerfeld's double covering of
three-space introduced in connection with the half-plane boundary.

Since one needs conformal invariance in two dimensions, not three, the \zfs\
are actually {\it modified} Barnes \zfs\ which have been dealt with in
[\pref{D2,Dow1,Cook}]. The determinants can be computed generally in terms
of Barnes \zfs\ but, because of the rational nature of $\pi/\be$, in this
case, they can be reduced to Epstein or Hurwitz \zfs. The general theory,
appropriate to the arbitrary 2-lune, is developed in [\pref{D2}]. However it
is probably easier to proceed directly.

From [\pref{D2}] the \zf\ for $-\De$ on the $ND$ 2-hemisphere is
 $$
\ze^{ND}_\pi(s)=\ze^{DD}_{2\pi}(s)-\ze^{DD}_\pi(s)
=\sum_{m,n=0}^\infty{1\over\big((1+m+n)^2-1/4\big)^s}\,.
 $$

Expanding in the $1/4$ leads to the expression for the derivative at 0,
 $$
{\ze^{ND}_{\pi}}'(0)=\ze_2'(0,1/2\mid1,1)+\ze_2'(0,3/2\mid1,1)
  -{N_2(1)\over4}\,,
\eql{deriv1}$$
where
 $$
\ze_2(s,a\mid1,1)=\sum_{m,n=0}^\infty{1\over(a+m+n)^s}
 $$
is a 2-dimensional Barnes \zf\ and $N_2(a)$ is its residue
at $s=2$; $N_2(a)=1$.

In this simple case the sums can easily be rearranged,
  $$
 \sum_{m,n=0}^\infty{1\over(a+m+n)^s}=\sum_{N=0}^\infty{N+1\over(N+a)^s}=
 \ze_R(s-1,a)+(1-a)\,\ze_R(s,a)\,,
 \eql{sums1}
  $$
so that
 $$
\ze_2(s,1/2\mid1,1)+\ze_2(s,3/2\mid1,1)=2\ze_R(s-1,1/2)=2(2^{s-1}-1)\,
\ze_R(s-1)
 $$
and therefore from (\peq{deriv1})
  $$
{\ze^{ND}_{\pi}}'(0)=-\ze_R'(-1)-{1\over12}\log2-{1\over4}.
  \eql{NDdet}
  $$
The absence of a $\ze_R'(0)$ term is related to the absence of the perimeter
\hk\ coefficient caused by the equal--sized $N$ and $D$ regions. The \zf\
has only the Weyl volume pole.

For comparison the standard formulae for the $DD$ and $NN$-hemispheres are
 $$
{\ze^{DD}_\pi}'(0)=2\ze_R'(-1)-\ze_R'(0)-{1\over4}
 $$
and
 $$
{\ze^{NN}_\pi}'(0)=2\ze_R'(-1)+\ze_R'(0)-{1\over4}.
 $$
By conformal transformation, on the $ND$-disc, the final result is
 $$
W_{ND}^{\rm disc}={1\over2}\ze_R'(-1)+{11\over24}\log2-{1\over24}\,,
 $$
using (\peq{cocycle2}) with (\peq{cocycle3}).
\section{\bf 7. The Robin boundary condition.\footnote{ According to
Gustafson, [\pref{Gustafson,GandA}], Gustave Robin (1855-1897), never seemed
to have used this condition. His name was, apparently, first attached to it
by Bergmann and Schiffer in the 1950's but the condition occurred already in
the work of Newton [\pref{newton}].}}

The Robin condition, (\peq{robin}), has made only a formal appearance in the
discussion so far. It was needed for conformal transformations but has not
yet entered into any eigenproblem.

The reason why Robin conditions are so awkward is that, in general, the
problem does not separate and, even if it does, the eigenvalues are given
only intrinsically. Early considerations of the eigenproblem are reviewed by
Pockels, [\pref{pockels}]. Poincar\'e, [\pref{Poincare}], also used the
condition in connection with eigenfunction existence. See also Bandle,
[\pref{bandle}]. A practical, more recent, treatment is given by Strauss,
[\pref{Strauss}].

Apart from applied mathematics, there has been some recent interest in Robin
conditions in the quantum field theoretic and spectral geometry scenes (\eg
Fulling, [\pref{fulling}], Bondurant and Fulling [\pref{BoandF}], Romeo and
Saharian [\pref{RandS}], Solodukhin, [\pref{Solod}], de Albuquerque and
Cavalcanti, [\pref{AandC}]).

In this section I return to the 2--hemisphere, some aspects of which were
mentioned in the previous section. Here I wish to see how far the standard
Robin eigenproblem on the interval, (\eg [\pref{pockels,Strauss,CandJ}]), is
relevant for a spherical geometry.

Appropriate details of the classic hemisphere eigenproblem  were also given
in [\pref{DGK}]. To repeat, cooordinates on the hemisphere are
$0\le\th\le\pi$ and $0\le\phi\le\pi$. The rim (boundary) consists of the
union of the two semicircles $\phi=0$ and $\phi=\pi$. In [\pref{DGK}], D
conditions were applied at $\phi=0$, and N at $\phi=\pi$. Now the latter are
replaced by Robin conditions and the former by either D or N. It is also
possible to treat both Robin, but I will not make this generalisation simply
for convenience. The singular region, $\Si$, comprises just the S and N
poles.

The Robin condition (\peq{robin}) specifically is
  $$
  {1\over\sin\th}{\pa\Phi\over\pa\phi}\bigg|_{\phi=\pi}
  =-S\Phi\bigg|_{\phi=\pi}\,,
  \eql{robin2}
  $$
which is not consistent with a separated structure for $\Phi$ unless $S$
takes the form
  $$
   S=-{h\over\sin\th}
   \eql{S}
  $$
diverging on $\Si$. I will return to this point later and take $S$ as in
(\peq{S}) simply in order to get on with the calculation because, in this
case, the condition (\peq{robin2}) reduces to the usual (D,R) (or (N,R)) on
the $\phi$ `interval' $(0,\pi)$ and I can employ known results. In the
separated solution for $\Phi$ the $\th$ part is unchanged, only the $\phi$
factor is modified. Thus, in the $(D,R)$ case, the hemisphere eigensolution
is
  $$
  \Phi_\la={\cal N}_h\sin (k\phi) P^{-k}_{n+k}(\cos\th)\,,
  \eql{hseigdr}
  $$
where $k>0$ is determined by the one--dimension interval condition,
  $$
k\cot(k\pi)=h\,,
  \eql{eigcond}
  $$
so that, as $h\to0$, $k$ tends to half an odd integer, $k\to m+1/2$,
$m=0,1,\ldots$, the Neumann result.

Likewise, for the $(N,R)$ case, instead of (\peq{hseigdr}) and
(\peq{eigcond}) there is
  $$
  \Phi_\la={\cal N}_h\cos (k\phi) P^{-k}_{n+k}(\cos\th)
  \eql{hseignr}
  $$
and
  $$
k\tan(k\pi)=-h\,.
  \eql{eigcondn}
  $$
This time, as $h\to0$, $k\to m$, $m=0,1,\ldots$. In both cases I will label
$k$ by the associated $m$.

To be specific, the eigenproblem I now wish to consider is
  $$
  (-\De+{1\over4})\Phi_\la=\la\Phi_\la\,,
  \eql{eigenp}
  $$
the reason being that the eigenvalues are perfect squares,
  $$
  \la=\la_{mn}=\big({1\over2}+k_m+n\big)^2\,,\quad m,n=0,1,\ldots\,.
  \eql{eigv2}
  $$
Degeneracies, if there are any, are due to coincidences.

For $(N,R)$ conditions with $h>0$, $k_0$ is imaginary but does not
correspond to a hemisphere mode because the assumed self--adjointness
implies that the eigenvalues $\la$ must be real. Then, in this particular
case, $m$ starts at 1.

One aim is to relate the spectral properties of the $\la$ to those of the
$k_m$ by summing out the $n$. I will do this via the heat and cylinder
kernels. We have used this before in spherical situations, [\pref{ChandD}].
The heat--kernel and cylinder (or `square root') kernel are defined, in
general, by
  $$\eqalign{
  K(t)&=\Tr e^{-Pt}=\sum_\la d_\la e^{-\la t}\cr
  K^{1/2}(t)&=\Tr e^{-\sqrt P\,t}=\sum_\la d_\la e^{-\sqrt\la\, t}\,,
  }
  \eql{kerns}
  $$
where I have included a degeneracy, just in case. From now on I use
Fulling's notation, setting $T(t)\equiv K^{1/2}(t)$, and taking $t$ as a
generic parameter.

Using the expression (\peq{eigv2}), it readily turns out that the hemisphere
cylinder kernels factorise,
  $$
  T_{HS}(t)= {1\over2\sinh t/2}\,T_I(t)\,,
  \eql{hsck}
  $$
$T_I$ being the cylinder kernel on the interval defined by,
  $$
  T_I(t)=\sum_{m=m_0}^\infty e^{-k_mt}\,,
  \eql{tid}
  $$
where $m_0=1$ for $(N,R)$ with $h>0$ and $m_0=0$ otherwise. Equation
(\peq{hsck}) is the connection between the hemisphere and the
$\pi$--interval.

Important information is contained in the short--time expansions of the heat
and cylinder kernels. A reflection of the pseudo-operator character of
$\big(-\De+1/4\big)^{1/2}$ is the possible existence of logarithmic terms in
the expansion of $T_{HS}$.

It is relatively straightforward to show that on a $d$--dimensional
manifold, maybe with a boundary, the expansion of a general cylinder kernel,
$T$, takes the form,
  $$
  T(t)\sim\sum_{i=0}^\infty a_i\, t^{-d+i}+\sum_{i}^\infty
  a'_i\,t^{-d+i}\log t
  \eql{expn1}
  $$

The lower limit on the second term is deliberately left unspecified  but I
draw attention to the important fact that, if the operator $P$ in
(\peq{kerns}) is a smooth differential operator with smooth boundary
conditions, then only {\it odd} positive powers of $t$ occur in the log
term, whatever the dimension of the manifold.

One way of showing this is to use the known existence, in this case, of the
series expansion of the heat--kernel $K(t)$, (\peq{hkasymp}), and then, by
\zf\ manipulations, relate the coefficients $a$, $a'$ and $C$. This was
first done in a physical context by Cognola {\it et al} [\pref{CVZ}]. A more
recent analysis, in the compact case, has been performed by B\"ar and
Moroianu [\pref{BandMo}] who consider local, diagonal kernels and give a
careful analysis of estimates. Another way is to relate the asymptotic
expansions obtained by smoothing using either $\la$ or $\sqrt\la$ as the
preferred variable. This was employed by Fulling [\pref{fulling}]. It is
sufficient to note, for now, that the coefficients of the logarithmic terms,
$a'_i$ in (\peq{expn1}) are determined by the $b_k$.

The heat--kernel coefficients, $b_k$, for $K_I$, on the $(R,R)$ interval
have been obtained in [\pref{Dow8}] and show that there are logarithmic
terms in $T_I$.  The relation (\peq{hsck}) between cylinder kernels then
implies that $T_{HS}$ also has logarithmic terms, but with {\it even} powers
of $t$. In turn, this suggests that the operator $-\De+1/4$, together with
the boundary conditions, including the choice of Robin parameter, (\peq{S}),
is a rather singular operator. This is looked at further in the next
section.
\section{\bf 8. Asymptotic series.}
To allow for the fact that the operator $-\De+1/4$ might be particularly
singular, I generalise (\peq{hkasymp}) to include logarithmic terms, the
immediate aim being to relate the expansions of the heat and cylinder
kernels. For this purpose, I use the zeta--function approach, mentioned
earlier, together with the general series established by Grubb and Seeley
[\pref{GrandS}]. Suitable summaries can be found in [\pref{Grubb2,Grubb3}].
For notational brevity I set $P=-\De+1/4$ and $Q=\sqrt P$, but, so far as
the general equations go, $P$ can be any Laplace--like (second order
elliptic) operator of smooth form.

The asymptotic expansion of the heat--kernel $K(t)$, (\peq{kerns}), is, \eg
[\pref{Grubb2,Grubb3}], on a $d$--dimensional manifold,
  $$
  K(t)\sim\sum_{-d\le k<-d+2} b_k\,t^{k/2}+\sum_{k=-d+2}^\infty\big(b'_k\log
  t+b''_k\big)\,t^{k/2}\,.
  $$
The reason for the lower limit of $-d+2$ will appear later. This limit
differs from that in [\pref{BandS,Grubb3}], which is zero.

The simple powers have been split into two sets because the coefficients
have different qualities. Since this does not concern me at this time, I
will combine them for algebraic ease. Therefore,
  $$
  K(t)\sim\sum_{k=-d}^\infty b_k\,t^{k/2}
  +\sum_{k=-d+2}^\infty b'_k\,t^{k/2}\log t\,.
  \eql{kasym}
  $$
This generalises (\peq{hkasymp}), with the relation between the coefficients
  $$
  b_k={C_{(k+d)/2}\over(4\pi)^{d/2}}\,.
  \eql{beec}
  $$
Similarly for the cylinder kernel ($Q$ is a first order pseudo--operator) I
will assume,
  $$
  T(t)\sim\sum_{k=-d}^\infty a_k\,t^k+\sum_{k=-d+2}^\infty a'_k\,t^k\log
  t\,,
  \eql{tasym}
  $$
where now, {\it all} powers of $t$ greater than $-d+1$ appear in the logarithmic
term.

The connection is made via the corresponding zeta--functions,
  $$
  \ze_Q(2s)=\ze_P(s)\,,
  \eql{reln}
  $$
which have asymptotic expansions corresponding to (\peq{kasym}) and
(\peq{tasym}). I refer to [\pref{Grubb2}], \eg, in order to save work, and
deduce,
  $$
  \Ga(s)\ze_P(s)\sim\sum_{k=-d}^\infty {b_k\over s+k/2}-{n_0\over s}-
  \sum_{k=-d+2}^\infty {b'_k\over(s+k/2)^2}
  \eql{zasymp}
  $$
  $$
  \Ga(s)\ze_Q(s)\sim\sum_{k=-d}^\infty {a_k\over s+k}-{n_0\over s}-
  \sum_{k=-d+2}^\infty {a'_k\over(s+k)^2}\,,
  \eql{zasymq}
  $$
where $n_0$ is the number of zero modes.

The relation between the $a$'s and the $b$'s follows from (\peq{reln}) using
the standard duplication formula for $\Ga(2s)$. Replacing $s$ by $2s$ in
(\peq{zasymq}) it is easy to see that dividing by $\Ga(s+1/2)$ removes
certain first order poles and converts some second order poles into first
order ones. From the residues, making the necessary identifications, I find
the relations
  $$\eqalign{
  b_k&={2^{k}\sqrt\pi\over\Ga\big((1-k)/2\big)}\,a_k\,,\quad k=-d,\ldots
  -1,0,2,4,\ldots\cr
  \noalign {\vskip5truept}
  &=(-1)^{(k+1)/2}\,2^{k-1}\,\Ga((k+1)/2)\,\sqrt\pi\,\,a'_k\,,\quad
  k=1,3,\ldots\cr
  \noalign {\vskip5truept}
  b'_k&={2^{k-1}\,\sqrt\pi\over\Ga\big((1-k)/2\big)}\,a'_k\,,\quad
  k=-d+2,\ldots,-1,0,2,\ldots\cr
  &=0\,\quad k=1,3,\ldots
  }
  \eql{bees}
  $$
which generalise those found by Fulling [\pref{fulling}] and Cognola {\it et
al}, [\pref{CVZ}].

Some overall conclusions can be drawn from these relations. Firstly, {\it
given} the heat--kernel $b_k$, it is not possible to determine the cylinder
$a_k$ for odd positive, $k=1,3,\ldots$, as emphasised by Fulling. Secondly,
for the assumed structure of the cylinder expansion (\peq{tasym}), or
(\peq{zasymq}), one cannot specify {\it all} the logarithmic terms in the
heat--kernel, \ie all the $b'_k$. Since the form (\peq{tasym}) is sufficient
for the quantities appearing in this paper, I will leave this point except
to say that it is easy to take (\peq{zasymp}) and work equation (\peq{reln})
the other way to derive the corresponding expansion for $\ze_Q(s)$. The new
feature is the appearance of poles of third order, leading to higher
logarithmic powers in the cylinder kernel expansion.

The extension of the lower limit down to $-d+2$ has the rather violent
consequence that, in this case, the coefficients $b_i$ for $i\ge -d+2$ are
global (\ie nonlocal) quantities. Only $b_{-d}$ and $b_{-d+1}$ are locally
computable. These are the Weyl volume and boundary area terms denoted by
$C_0$ and $C_{1/2}$ earlier. This means, in particular that $C_1$, discussed
in \S3, is not locally computable. We have seen that this is true for the
exhibited form, (\peq{cee1}), through the last two, \ie codimension 2,
terms. The question of nonlocal terms is considered in a later section.

A further consequence of (\peq{bees}), and the one relevant for the
spherical problem treated in the previous section, is that the logarithmic
terms in the heat--kernel arise from those logarithmic terms in the cylinder
kernel with {\it even} powers of $t$, which is precisely the case with
(\peq{hsck}).

The conclusion is that the Laplace operator, $-\De+1/4$ with the boundary
conditions, {\it is} a rather singular operator as it provides a concrete
example of a second order problem involving logarithmic terms. It would
therefore seem that it does not come within the compass of Seeley's
analysis, [\pref{Seeley}], probably because of the divergence in the Robin
function, $S$. Hence, before giving explicit forms for the expansions, I
break off to look at the modes, (\peq{hseigdr}), a little more closely, one
reason being that, although self--adjointness depends on the formal
subtraction of two (identical) terms,
  $$
  \int_0^\pi d\th\,\Phi\,S\,\Psi\,,
  \eql{bint}
  $$
at the boundary, the divergence of $S$, (\peq{S}), at the poles
($\th=0,\pi$), \ie on $\Si$, might give one pause for thought. In general
terms, the assumed Dirichlet conditions at $\Si$ are actually sufficient for
convergence of the integral, (\peq{bint}).
\section{\bf 9. Robin mode properties.}
The orthonormality of (\peq{hseigdr}), and of (\peq{hseignr}), is easily
established either by direct integration or, formally, by the usual
self--adjoint Liouville method.

For convenience I write down some standard things. For the Legendre
functions orthonormality reads
  $$
  \int_{-1}^1 dx\,P^{-k}_{n+k}(x)\,P^{-k}_{n'+k}(x)={2\Ga(n+1)\over
 (2k+2n+1)  \Ga(2k+n+1)}\,\de_{nn'}\,,
  $$
where $n$ and $n'$ are positive integers or zero and $k>-1 $ (MacRobert
[\pref{MacRobert}] p.335).

I also note the limiting behaviour at the poles,
  $$
  P^{-k}_{n+k}(z)\to (1-z^2)^{k/2}{1\over 2^k\Ga(k+1)}\,,\quad z\to\pm1\,.
  \eql{lim1}
  $$

The interval Robin modes are standard, (\eg Carslaw and Jaeger,
[\pref{CandJ}], Strauss [\pref{Strauss}], Pockels, [\pref{pockels}]) with
easily determined normalisations. I again remark that Dirichlet
conditions continue to apply on $\Si$, my default position.

Use of the limit (\peq{lim1}) shows that each boundary term, (\peq{bint}),
in the
self--adjoint condition applied to two eigenfunctions, corresponding to $k$
and $k'$, is finite if $k+k'>0$. So we are completely safe in this case.

This limit also implies that the total heat flux, per mode,
  $$
  ~\,\,\int_{\pa\man}\pa_n\Phi_\la\,,
  $$
is finite. Moreover, the `quantum mechanical flux',
  $$
   ~\,\,\int_{\pa\man}\Phi_\la\pa_n\Phi_\la\,,
  $$
is also finite, by virtue of Barnes' formula,
  $$
  \int_0^1dx\,{\big(P^\mu_\nu(x)\big)^2\over1-x^2}=
  -{1\over2\mu}{\Ga(1+\mu+\nu)\over\Ga(1-\mu+\nu)}\,,
  \eql{barnes}
  $$
valid for $\Real \mu<0$ and $\mu+\nu$ a positive integer, or zero. The lower
limit can be extended to $-1$ using the fact that the Legendre functions are
unchanged, up to a sign, under $x\to-x$ ([\pref{MacRobert}], p.334 Examples
(2) and (3)).

Confidence in the eigenmodes is increased if one employs a perturbation
technique to calculate the change in the eigenvalues $\la_{mn}$ of
(\peq{eigv2}) for a small change in $S$. The formula, which will not be
derived here, is
  $$
  \de\la=\int_{\pa\man(N)}\Phi_\la S\Phi_\la\,.
  \eql{pert}
  $$
For simplicity, I am considering the $(D,R)$ set--up and am perturbing about
the $(D,N)$ case, so $S$ is small. The integration over the boundary
encompasses only the part on which $R$ (equivalently $(N)$) holds, since
$\Phi$ is zero on the $D$ part. Equation (\peq{pert}) seems to occur first,
for constant $S$, in Poincar\'e [\pref{Poincare}] and for variable $S$ in
Pockels, [\pref{pockels}] p.178, as a quick consequence of Green's formula.
See also Fr\"ohlich, [\pref{Frohlich}], equn. (6d).

The answer is known from a direct analysis of the interval eigenvalue
condition (\peq{eigcond}) which shows that
   $$
   k_m-{1\over2}\sim m-{2h\over(2m+1)\pi}\,.
   \eql{change}
   $$

Alternatively, applying (\peq{pert}) yields
  $$
  \de\la_{mn}=-h\can^2\int_0^\pi d\th\,{1\over\sin\th}\big(P^{-\ol m/2}_{n+\ol
  m/2}(\cos\th)\big)^2\,,
  \eql{pert2}
  $$
where $\ol m =2m+1$. The normalisation is,
  $$
  \can^2={2\over\pi}{\Ga(\ol m+n+1)\Ga(\ol m+2n+1)\over2\Ga(n+1)}\,,
  $$
and substitution of (\peq{barnes}) into (\peq{pert2}) easily produces
  $$
   \de\sqrt\la_{mn}=-{2h\over\pi}{1\over\ol m}\,,
  $$
which is just (\peq{change}). This limited check inspires a certain
confidence in the sensibility of the modes and the model.

I should point out now that I do not attach any serious significance to the
model. The choice $S$ is simply one of convenience for solvability.
\section{\bf 10. Explicit expressions}

In this section I use the form of the heat--kernel coefficients for the
interval $(R,R)$ Robin problem derived in [\pref{Dow8}], from which one can
easily deduce those for the $(D,R)$ and $(N,R)$ cases by appropriate limits.

I find, for both $(N,R)$ ($h<0$) and $(D,R)$ on an interval of length $\pi$,
(I), the heat--kernel coefficients,
  $$
  b_k=b^I_k={1\over2\Ga(k/2+1)}\,h^k\,,\quad k=1,\ldots,
  \eql{rbees}
  $$
so that, for example, the coefficients of the logarithmic terms in the
interval {\it cylinder} kernel, (\peq{tid}), are determined from
(\peq{bees}) applied to the interval, as
  $$
  a'_{2n-1}=(-1)^n{1\over\pi}\,{ 1
  \over(2n-1)!}\,h^{2n-1}\,,\quad n=1,2,\ldots\,,
  $$
and these odd terms are all there are because there are no logarithmic terms
in the interval heat--kernel. The asymptotic series can be summed to produce
the closed form for the asymptotic logarithmic part of $T_I$,
  $$
  T^{log}_I(t)\sim {e^{-h^2t^2}-1\over \pi ht}\,\log t\,,
  $$
and then from (\peq{hsck}) the logarithmic part of the $(N,R)$ ($h<0$), or
$(D,R)$, hemisphere, (HS), cylinder kernel follows as,
  $$
  T^{log}_{HS}(t)\sim \,{1\over\sinh(t/2)}
  {e^{-h^2t^2}-1\over 2\pi ht}\,\log t\,.
  \eql{tlog}
  $$

The same result also holds for the $(N,R)$, $h>0$, case even after dropping
the imaginary $k_0$ interval--mode. The contribution of this to the \hk\
coefficients must be removed but this does not affect the $b^I_k$, for odd
$k$.

The relation (\peq{bees}) can now be applied to the hemisphere and the
logarithmic terms in the {\it heat--kernel} worked out. I do not give the
easily derived expressions. They are combinations of Bernoulli numbers.
Their existence is all that is required for now. I only remark that the
important term $\sim t^0\log t$ is present.

There are, of course, a number of technical routes to the expressions and
conclusions derived above. I have chosen to use heat and cylinder kernels
but one could employ the resolvent and a standard contour method of
rewriting eigenvalue sums. I also note that it is straightforward to extend
the calculations to a lune and also to $d$--dimensions, which just gives
higher powers of $\sinh(t/2)$ in (\peq{tlog}), say. This justifies the lower
limit in (\peq{tasym}) and hence in (\peq{kasym}). Odd spheres are not
singular, in agreement with a general result. The details will be presented at
another time.
\section{\bf 11. Nonlocal terms and the Casimir energy}
In this section I enlarge on previous statements regarding the nonlocality
of some expansion coefficients, in particular of $C_1$ which, in two
dimensions, is the coefficient of the constant term in the \hk\ expansion.

As a first step, I look at the Casimir energy, $E$, on the $\pi$--interval
(I) and simply quote the formal definition,
  $$
  E={\rm FP}\,{1\over2}\,\ze^I_P(-1/2)\,,
  \eql{cas}
  $$
where $P=-d^2/d\phi^2$ and the boundary conditions are either (D,R) or (N,R)
on the ends. FP stands for `finite part'. The Casimir energy is a nonlocal
quantity.

On the $\pi$--interval the relevant definition here of the \zf\ is,
  $$
  \ze^I_P(s)=\sum_{m=m_0}{1\over k_{m}^{2s}}\,.
  \eql{zef}
  $$

For the purposes of this paper, I consider $E$, defined by (\peq{cas}),
simply as a convenient mathematical quantity rather than as something having
physical, operational significance. In particular, on the pure interval,
there would be the question of what to do with the `negative mode' that
occurs for $(N,R)$ when $h>0$. In the present model, this has already been
excluded. (Compare with reference [\pref{Solod}] where these modes are
referred to as `bound states'.)

Equations (\peq{reln}), (\peq{zasymp}) and (\peq{zasymq}) hold for the
interval. Since $P$ is smooth, there are no log terms in the \hk\ and so the
coefficients $b^{'I}_k$ are zero.

It is important to realise that in the Robin case, $\ze_P^I(s)$ has a pole
at $s=-1/2$,
  $$
  \ze_P^I(s)\sim  {A\over s+1/2}+B
  \eql{ize}
  $$
where the residue, $A=-{b_1^I/2\sqrt\pi}$, follows from (\peq{zasymp}) and
the remainder, $B$, equals $2E$, by definition. In terms of the interval
\hk\ coefficients, $C^I_n$, (see (\peq{rbees}), (\peq{beec})),
  $$
  b_1^I={1\over2\sqrt\pi}\,C^I_1={h\over\sqrt\pi}\,,
  \eql{bee1}
  $$
for both $(D,R)$ and $(N,R)$, for all $h$.

From (\peq{cas}) and (\peq{reln}) it is required to work around $s\sim-1$
for $\ze_Q^I(s)$ where
  $$
  \Ga(s)\ze_Q^I(s)\sim {a_1^I\over s+1}-{a_1^{'I}\over(s+1)^2}\,.
  \eql{zeq}
  $$
and so
  $$
  \bigg(-{1\over s+1}+\ga-1\bigg)\bigg({A\over s/2+1/2}+B\bigg)
  \sim {a_1^I\over s+1}-{a_1^{'I}\over(s+1)^2}\,,
  $$
which yields
  $$
  b_1^I=-\sqrt\pi\,a_1^{'I}\,,
  \eql{ba}
  $$
and also
  $$
  E=-{1\over2}\,a_1^I+{\ga-1\over2\sqrt\pi}\,b_1^I
  =-{1\over2}\,a_1^I-{\ga-1\over2\pi}\,h\,.
  \eql{cas2}
  $$
Equation (\peq{ba}) agrees with (\peq{bees}) while (\peq{cas2}) relates the
coefficient $a_1^I$ to the nonlocal Casimir energy. It differs from the
result in [\pref{fulling}] by the last (constant) term. See [\pref{CVZ}] for
a relevant discussion of various regularisation recipes in this context.

The next step is to relate the interval and hemisphere expansions using
(\peq{hsck}). Expansion of the $1/\sinh$ easily gives the connection,
  $$
  a_0^{HS}=a_1^I-{1\over24}a_{-1}^I\,,
  \eql{conn1}
  $$
which I now rewrite in terms of \hk\ expansion coefficients.

Application of (\peq{bees}) to the two--hemisphere gives,
  $$
  a_0^{HS}=b_0^{HS}\,,
  $$
where $b_0^{HS}$ is the constant term coefficient,
  $$
  b_0^{HS}={1\over4\pi}\,C_1^{HS}\,.
  $$

Equation (\peq{bees}), applied to the interval, gives,
  $$
  a_{-1}^I={2\over\sqrt\pi}\,b_{-1}^I
  ={2\over\sqrt\pi}{C_0^I\over2\sqrt\pi} =1\,,
  $$
and so, finally, (\peq{conn1}) becomes, using (\peq{cas2}),
  $$
  C_1^{HS}=-8\pi\,E-{\pi\over6}-4h\,(\ga-1)\,,
  \eql{conn2}
  $$
which relates the $C_1$ coefficient on the hemisphere to the Casimir energy
on the interval {\it for my model}, \ie (\peq{S}).

A basic check sets $h=0$ when (\peq{conn2}) allows a computation of the
interval Casimir energies for the (D,D) and (D,N) cases using the
expressions for $C_1$, (\peq{cee1d}), (\peq{cee1r}) and (\peq{cee1}).
Simple algebra yields the standard values,
  $$
  E(D,D)=E(N,N)=-{1\over24}\,,\quad E(D,N)={1\over48}\,.
  $$

It is worth noting that for the (D,D) and (N,N) cases, although the Casimir
energy, $E$, is nonlocal on the interval, the $C_1$ coefficients are
local on the hemisphere.

In the general case, $E$ is a nontrivial function of $h$ on the interval and
the conjectured forms of $C_1(D,R)$ and $C_1(N,R)$, (\peq{cee1}),
(\peq{cee1n}), obviously do not agree with (\peq{conn2}). However, it will
be recalled that the \hk\ expansion is really one in $h^2 t$ and therefore
one should treat $h$ as `small'. As mentioned, it is in this realm that
(\peq{cee1}) and (\peq{cee1r}) should be valid. Furthermore the divergence
of the hemisphere Robin function, $S$ of (\peq{S}), makes $C_1$ formally
infinite, and cannot be considered `small'. Hence it is clear that
(\peq{cee1}), (\peq{cee1r}) and (\peq{cee1n}) are not even appropriate for
the present situation. Against this must be set the fact that, as described
in section 9, perturbation in $S$ appears to work for the eigenvalues.
Therefore I intend to give a further look at perturbation theory. This will
give us a simple, if limited, analytical handle on the \zf\ on the Robin
interval which might also be useful in other circumstances.

\section {\bf 12. Perturbation approach}

I first consider the (N,R) case when one must distinguish between positive
and negative Robin parameters, $h$. If $h>0$  one has the approximation for
the interval frequencies,
  $$
  k_m\approx m-{h\over m\pi}\,,\quad m=1,2,\ldots\,,
  $$
which leads to
  $$
  \ze^I_P(s)\approx\ze_R(2s)+{2hs\over\pi}\,\ze_R(2s+2)\,,
  \eql{zenr}
  $$
in terms of the Riemann \zf, and so $A=-h/2\pi$ with
  $$
  E(N,R)\approx-{1\over24}-{h\over2\pi}\,(\ga-1)\,,\quad h\downarrow0\,.
  \eql{enr1}
  $$

By contrast, if $h<0$, there is a root which tends to zero as $h\to0$ and
becomes the zero (N,N) mode. Approximation of (\peq{eigcondn}) gives,
  $$\eqalign{
      k_0\approx&\, \sqrt{-h\over\pi}\,,\quad h\uparrow0\cr
      k_m\approx&\, m-{h\over m\pi}\,,\quad m=1,2,\ldots\,,
      }
  $$
whence $A=-h/2\pi$ and
  $$
  E(N,R)\approx -{1\over24}+{1\over2}\sqrt{-h\over\pi}
  -{h\over2\pi}\,(\ga-1)\,,\quad h\uparrow0\,.
  \eql{enr2}
  $$

Turning to the (D,R) case, to order $h$, (\peq{change}),
  $$
  k_m\approx{\ol m\over2}-{2h\over\ol m\,\pi}\,
  \quad \ie \quad \la\approx{\ol m^2\over4}-{2h\over\pi}\,,\quad \ol
  m=1,3,\ldots\,,
  $$
so that
  $$
  \ze^I_P(s)\approx (2^{2s}-1)\,\ze_R(2s)+{2hs\over\pi}\,
  (2^{2s+2}-1)\ze_R(2s+2)\,.
  \eql{zedr}
  $$
Working around $s\sim-1/2$ yields, in accordance with (\peq{ize}), the
residue check, $ A=-h/2\pi$, and the Casimir energy,
  $$
  E(D,R)\equiv{B\over2}\approx{1\over48}-{h\over2\pi}\big(
  \ga-1+2\log2\big)\,.
  \eql{cas3}
  $$

From these expressions, I use (\peq{conn2}) to {\it compute} the \hk\
coefficients on the hemisphere. I find, in the two cases, the values,
  $$\eqalign{
  C_1(N,R)&\approx{\pi\over6}\,,\quad h>0\cr
  C_1(N,R)&\approx{\pi\over6}-4\sqrt{-\pi h}\,,\quad h<0\cr
  C_1(D,R)&\approx-{\pi\over3}+8h\log2\,.
  }
  \eql{cee1h}
  $$
It appears that the effect of the Robin function (\peq{S}) for $h>0$ has
disappeared from $C_1(N,R)$ leaving just the $C_1(N,N)$ value, at least in
this lowest approximation. This asymmetry can be traced to the omission of
the negative mode for $h>0$. (Compare with the discussion in
[\pref{Solod}].)

The log 2 term in $C_1(D,R)$ reinforces the conclusion that (\peq{cee1}) is
inappropriate for the present singular model.

Similar results can be shown to hold for the simple wedge of section 2 with
Robin conditions on one side. In order for the techniques of [\pref{BKD}] to
work, separability demands that the Robin function is again singular on
$\Si$ (the apex of the wedge) taking the form $S=-h/r$. Without going into
details, just referring to equations (4.4)-(4.6) of [\pref{BKD}], there is
again a log term coming from the pole at $s=-1/2$ in the interval Robin \zf\
and there is a log 2 term in the expression for the corresponding $C_1$.

\section {\bf 13. Exact form of interval Casimir energy.}

The exact expressions for the Casimir energy derived in [\pref{RandS}] and
[\pref{AandC}] can be approximated for comparison with my perturbation
results. For convenience I refer to equations (34) and (39) in
[\pref{AandC}] for the Casimir energies and rewrite them in one dimension,
   $$\eqalign{
   \caE(N,R)=E(N,N)+{1\over2\pi}\int dk
   \log\bigg(1-{2h\over(k-h)(\exp(2k\pi-1))}\bigg)\cr
   \caE(D,R)=E(D,N)+{1\over2\pi}\int dk
   \log\bigg(1+{2k\over(k-h)(\exp(2k\pi-1))}\bigg)-{1\over16}\,,
   }
   \eql{intcasen}
   $$
where the $\caE$ may, or may not, agree, up to renormalisation, with the
quantities evaluated from (\peq{cas}). (I have set $h=-c_2<0$, in the
notation of [\pref{AandC}], to comply with my sign.)

The leading small $h$ behaviours of the `correction' terms in
(\peq{intcasen}) can be determined {\it numerically} to be,
  $$\eqalign{
   \caE(N,R)&\approx E(N,N)+ {1\over 2}\,\sqrt {-h\over\pi}
   +{h\over2\pi}\,F_N(h)\cr
   \caE(D,R)&\approx E(D,N)+{h\over2\pi}\,F_D(h)\,,
   }
   \eql{intcasen2}
   $$
where $F_N$ and $F_D$ satisfy the functional relation
  $$
  F(\la x)-\la F(x)\approx{\log\la\over\la}\,.
  \eql{fr}
  $$
There is some similarity between (\peq{intcasen2}) and the perturbative
results, but I cannot explain (\peq{fr}) from such a viewpoint and simply
present the result as a possible significant curiosity.

\begin{ignore}

I now turn to the (D,R)--{\it hemisphere} \zf\ constructed from the
eigenvalues, (\peq{eigv2}),
  $$\eqalign{
  \ze^{HS}(s)&=\sum_{m,n=0}^\infty{1\over(1/2+k_m+n)^{2s}}\cr
  &=\ze_R(2s-1)+\de\ze^{HS}(s)
}
  $$
where (\peq{sums1}) has been used for the unperturbed (D,N)--hemisphere \zf.
The small change in the \zf\ is,
  $$\eqalign{
  \de\ze^{HS}(s)&=-{4sh\over\pi}\sum_{m,n=0}^\infty {1\over 2m+1}{1\over(
  1+m+n)^{2s+1}}\cr
  &=-{4sh\over\pi}\sum_{N=0}^\infty \sum_{m=0}^N{1\over 2m+1}{1\over(
  1+N)^{2s+1}}\,.
}
  $$

The aim is to evaluate $C_1$ from the standard relation, in two dimensions,
  $$
  \ze(0)={C_1\over4\pi}\,,
  $$
in the absence of zero modes. It should be a formality to regain the value
in (\peq{cee1h}) since the analysis is only being done in a different order.
\end{ignore}
\section{\bf 14. Conclusion}

Apart from rectilinear domains, and the hemisphere, there seem few
situations that can be solved exactly for $ND$-conditions (see
[\pref{pockels}]) and this is a drawback to the construction of the precise
forms of the \hk\ coefficients. Nevertheless a certain progress has been
made in a simple minded way making use of the $ND$-wedge expression. This
type of reasoning can be extended to higher dimensions leading to
information about the trihedral corner contributions and their higher
analogues.

Surprisingly the conformal functional determinant is available on the
`half-N half-D' disc by conformal transformation from the $ND$-hemisphere,
and has been computed, assuming that a conjecture for the heat--kernel
coefficient, $C_1$, is correct.

In sections 7 to 10, I considered hybrid $(D,R)$ and $(N,R)$ 2--hemisphere
problems, with a Laplace style operator, which have logarithmic terms in the
short--time expansions of the heat--kernels. The boundary conditions involve
a Robin function that diverges, but not too strongly, at the poles which are
the places where the $D,N$ conditions meet the $R$ condition. This allows
the problem to be separated.

Although the second order differential operator (with the boundary
conditions) is rather singular, it is not singular enough to prevent the
heat--kernel expansion from existing.

The (singular) Robin condition is essential here for the existence of the
logarithmic terms. The expression (\peq{tlog}) vanishes when $h=0$. In
[\pref{DGK}] it was suggested, by an indirect argument, that for an $(N,D)$
situation, a nonzero extrinsic curvature at the boundary gave rise to
logarithmic terms. In general, it seems that the Robin condition mimics an
extrinsic curvature, [\pref{fulling}] equn, (4.1). If this is so, it is
perhaps not surprising that a singular Robin function, $S$, produces a log
term as it would simulate a conical type singularity.

The attempt to relate the hemisphere \hk\ coefficient to the Robin Casimir
energy on the interval, analysed in sections 11 and 12, confirms the limited
validity of the expression for the $C_1$ coefficient, (\peq{cee1}),
to extend which requires further analysis.

An incidental line of technical enquiry would centre on the significance of
the `imaginary modes' on the $(N,R)$ interval for $h>0$ and whether
continuity could be maintained through $h=0$.

Further aspects of this particular set--up will be examined at another time.

\section{\bf References}
\vskip 5truept
\begin{putreferences}
 \ref{AandD}{Apps,J.S. and Dowker,J.S. \cqg{15}{1998}{1121}.}
 \ref{APS}{Atiyah,M.F., V.K.Patodi and I.M.Singer: Spectral asymmetry and
Riemannian geometry \mpcps{77}{1975}{43}.}
 \ref{avramidi}{Avramidi,I. {\it
Heat kernel asymptotics of a non-smooth boundary value problem} Workshop on
Spectral Geometry, Bristol 2000.}
 \ref{AandT}{Awada,M.A. and D.J.Toms:
Induced gravitational and gauge-field actions from quantised matter fields
in non-abelian Kaluza-Klein thory \np{245}{1984}{161}.}
 \ref{BandI}{Baacke,J.
and Y.Igarishi: Casimir energy of confined massive quarks
\prD{27}{83}{460}.}
 \ref{Barnesa}{Barnes,E.W.: On the Theory of the
multiple Gamma function {\it Trans. Camb. Phil. Soc.} {\bf 19} (1903) 374.}
\ref{Barnesb}{Barnes,E.W.: On the asymptotic expansion of integral
functions of multiple linear sequence, {\it Trans. Camb. Phil. Soc.} {\bf
19} (1903) 426.}
 \ref{Barv}{Barvinsky,A.O. Yu.A.Kamenshchik and
I.P.Karmazin: One-loop quantum cosmology \aop {219}{1992}{201}.}
 \ref{BandM}{Beers,B.L. and Millman, R.S. :The spectra of the
Laplace-Beltrami operator on compact, semisimple Lie groups.
\ajm{99}{1975}{801-807}.}
 \ref{BandH}{Bender,C.M. and P.Hays: Zero point
energy of fields in a confined volume \prD{14}{1976}{2622}.}
 \ref{BBG}{Bla\v
zi\' c,N., Bokan,N. and Gilkey,P.B.: Spectral geometry of the form valued
Laplacian for manifolds with boundary \ijpam{23}{1992}{103-120}}
 \ref{BEK}{Bordag,M., E.Elizalde and K.Kirsten: { Heat kernel coefficients
of the Laplace operator on the D-dimensional ball}, \jmp{37}{1996}{895}.}
 \ref{BGKE}{Bordag,M., B.Geyer, K.Kirsten and E.Elizalde,: { Zeta function
determinant of the Laplace operator on the D-dimensional ball},
\cmp{179}{96}{215}.}
 \ref{BKD}{Bordag,M., K.Kirsten,K. and Dowker,J.S. \cmp{182}{96}{371}.}
 \ref{BKD1}{Bordag,M., K.Kirsten,K. and Dowker,J.S.: Heat
kernels and functional determinants on the generalized cone
\cmp{182}{96}{371}.} \ref{Branson}{Branson,T.P.: Conformally covariant
equations on differential forms \cpde{7}{1982}{393-431}.}
 \ref{BandG2}{Branson,T.P. and Gilkey,P.B. {\it Comm. Partial Diff. Eqns.}
{\bf 15} (1990) 245.}
 \ref{BandG}{Branson,T.P. and Gilkey,P.B.
\tams{344}{1994}{479}.}
 \ref{BGV}{Branson,T.P., P.B.Gilkey and
D.V.Vassilevich {\it The Asymptotics of the Laplacian on a manifold with
boundary} II, hep-th/9504029.}
 \ref{BCZ1}{Bytsenko,A.A, Cognola,G. and
Zerbini, S. : Quantum fields in hyperbolic space-times with finite spatial
volume, hep-th/9605209.}
 \ref{BCZ2}{Bytsenko,A.A, Cognola,G. and Zerbini,
S. : Determinant of Laplacian on a non-compact 3-dimensional hyperbolic
manifold with finite volume, hep-th /9608089.}
 \ref{CandH2}{Camporesi,R.
and Higuchi, A.: Plancherel measure for $p$-forms in real hyperbolic space,
\jgp{15}{1994}{57-94}.}
 \ref{CandH}{Camporesi,R. and A.Higuchi {\it On the
eigenfunctions of the Dirac operator on spheres and real hyperbolic spaces},
gr-qc/9505009.} \ref{ChandD}{Chang, Peter and Dowker,J.S.
\np{395}{1993}{407}.} \ref{cheeg1}{Cheeger, J. \jdg {18}{1983}{575}.}
\ref{cheeg2}{Cheeger,J.: Hodge theory of complex cones {\it Ast\'erisque}
{\bf 101-102}(1983) 118-134} \ref{Chou}{Chou,A.W.: The Dirac operator on
spaces with conical singularities and positive scalar curvature,
\tams{289}{1985}{1-40}.} \ref{CandT}{Copeland,E. and Toms,D.J.: Quantized
antisymmetric tensor fields and self-consistent dimensional reduction in
higher-dimensional spacetimes, \break\np{255}{1985}{201}} \ref{Cook}{Cook,A.
1986 {\it PhD Thesis}, University of Manchester.} \ref{DandH}{D'Eath,P.D.
and J.J.Halliwell: Fermions in quantum cosmology \prD{35}{1987}{1100}.}
\ref{cheeg3}{Cheeger,J.:Analytic torsion and the heat equation. \aom{109}
{1979}{259-322}.} \ref{DandE}{D'Eath,P.D. and G.V.M.Esposito: Local boundary
conditions for Dirac operator and one-loop quantum cosmology
\prD{43}{1991}{3234}.} \ref{DandE2}{D'Eath,P.D. and G.V.M.Esposito: Spectral
boundary conditions in one-loop quantum cosmology \prD{44}{1991}{1713}.}
\ref{D}{Dowker,J.S. \cqg{11}{1994}{557}.} \ref{Dow1}{Dowker,J.S.
\cmp{162}{1994} {633}.}
 \ref{Dow8}{Dowker,J.S. \cqg{13}{1996}{585}; hep-th/9506042.}
 \ref{Dow9}{Dowker,J.S. {\it Oddball determinants} MUTP/95/12;
hep-th/9507096.} \ref{Dow10}{Dowker,J.S. {\it Spin on the 4-ball},
hep-th/9508082, {\it Phys. Lett. B}, to be published.}
\ref{DandA}{Dowker,J.S. and Apps,J.A. \cqg{12}{1995}{1363}.}
\ref{D2}{Dowker,J.S. \jmp{35}{1994}{4989}; 1995 (Feb.) erratum.}
\ref{DandA2}{Dowker,J.S. and J.S.Apps, {\it Functional determinants on
certain domains}. To appear in the Proceedings of the 6th Moscow Quantum
Gravity Seminar held in Moscow, June 1995; hep-th/9506204.}
\ref{DABK}{Dowker,J.S., Apps,J.S., Bordag,M. and Kirsten,K.: Spectral
invariants for the Dirac equation with various boundary conditions {\it
Class. Quant.Grav.} to be published, hep-th/9511060.}
\ref{EandR}{E.Elizalde
and A.Romeo : An integral involving the generalized zeta function, {\it
International J. of Math. and Phys.} {\bf13} (1994) 453.}
\ref{ELV2}{Elizalde, E., Lygren, M. and Vassilevich, D.V. : Zeta function
for the laplace operator acting on forms in a ball with gauge boundary
conditions. hep-th/9605026}
\ref{ELV1}{Elizalde, E., Lygren, M. and
Vassilevich, D.V. : Antisymmetric tensor fields on spheres: functional
determinants and non-local counterterms, \jmp{}{96}{} to appear. hep-th/
9602113} \ref{Kam2}{Esposito,G., A.Y.Kamenshchik, I.V.Mishakov and
G.Pollifrone: Gravitons in one-loop quantum cosmology \prD{50}{1994}{6329};
\prD{52}{95}{3457}.} \ref{Erdelyi}{A.Erdelyi,W.Magnus,F.Oberhettinger and
F.G.Tricomi {\it Higher Transcendental Functions} Vol.I McGraw-Hill, New
York, 1953.} \ref{Esposito}{Esposito,G.: { Quantum Gravity, Quantum
Cosmology and Lorentzian Geometries}, Lecture Notes in Physics, Monographs,
Vol. m12, Springer-Verlag, Berlin 1994.}
\ref{Esposito2}{Esposito,G. {\it
Nonlocal properties in Euclidean Quantum Gravity}. To appear in Proceedings
of 3rd. Workshop on Quantum Field Theory under the Influence of External
Conditions, Leipzig, September 1995; gr-qc/9508056.}
\ref{EKMP}{Esposito G,
Kamenshchik Yu A, Mishakov I V and Pollifrone G.: One-loop Amplitudes in
Euclidean quantum gravity. \prd {52}{1996}{3457}.}
\ref{ETP}{Esposito,G.,
H.A.Morales-T\'ecotl and L.O.Pimentel {\it Essential self-adjointness in
one-loop quantum cosmology}, gr-qc/9510020.}
\ref{FORW}{Forgacs,P.,
L.O'Raifeartaigh and A.Wipf: Scattering theory, U(1) anomaly and index
theorems for compact and non-compact manifolds \np{293}{1987}{559}.}
\ref{GandM}{Gallot S. and Meyer,D. : Op\'erateur de coubure et Laplacian des
formes diff\'eren-\break tielles d'une vari\'et\'e riemannienne
\jmpa{54}{1975} {289}.} \ref{Gilkey1}{Gilkey, P.B, Invariance theory, the
heat equation and the Atiyah-Singer index theorem, 2nd. Edn., CRC Press,
Boca Raton 1995.} \ref{Gilkey2}{Gilkey,P.B.:On the index of geometric
operators for Riemannian manifolds with boundary \aim{102}{1993}{129}.}
\ref{Gilkey3}{Gilkey,P.B.: The boundary integrand in the formula for the
signature and Euler characteristic of a manifold with boundary
\aim{15}{1975}{334}.} \ref{Gottlieb}{Gottlieb, H.P.W., {\it
J.Aust.Math.Ass.} {\bf26} (1984) 293.} \ref{Grubb}{Grubb,G. {\it Comm.
Partial Diff. Eqns.} {\bf 17} (1992) 2031.} \ref{GandS1}{Grubb,G. and
R.T.Seeley \cras{317}{1993}{1124}; \invm{121}{1995} {481}.}
\ref{GandS}{G\"unther,P. and Schimming,R.:Curvature and spectrum of compact
Riemannian manifolds, \jdg{12}{1977}{599-618}.} \ref{IandT}{Ikeda,A. and
Taniguchi,Y.:Spectra and eigenforms of the Laplacian on $S^n$ and $P^n(C)$.
\ojm{15}{1978}{515-546}.} \ref{IandK}{Iwasaki,I. and Katase,K. :On the
spectra of Laplace operator on $\La^*(S^n)$ \pja{55}{1979}{141}.}
\ref{JandK}{Jaroszewicz,T. and P.S.Kurzepa: Polyakov spin factors and
Laplacians on homogeneous spaces \aop{213}{1992}{135}.}
\ref{Kam}{Kamenshchik,Yu.A. and I.V.Mishakov: Fermions in one-loop quantum
cosmology \prD{47}{1993}{1380}.} \ref{KandM}{Kamenshchik,Yu.A. and
I.V.Mishakov: Zeta function technique for quantum cosmology {\it Int. J.
Mod. Phys.} {\bf A7} (1992) 3265.} \ref{KandC}{Kirsten,K. and Cognola.G,: {
Heat-kernel coefficients and functional determinants for higher spin fields
on the ball} \cqg{13}{1996} {633-644}.} \ref{Levitin}{Levitin,M.: {
Dirichlet and Neumann invariants for Euclidean balls}, {\it Diff. Geom. and
its Appl.}, to be published.} \ref{Luck}{Luckock,H.C.: Mixed boundary
conditions in quantum field theory \jmp{32}{1991}{1755}.}
\ref{MandL}{Luckock,H.C. and Moss,I.G,: The quantum geometry of random
surfaces and spinning strings \cqg{6}{89}{1993}.} \ref{Ma}{Ma,Z.Q.: Axial
anomaly and index theorem for a two-dimensional disc with boundary
\jpa{19}{1986}{L317}.} \ref{Mcav}{McAvity,D.M.: Heat-kernel asymptotics for
mixed boundary conditions \cqg{9}{1992}{1983}.} \ref{MandV}{Marachevsky,V.N.
and D.V.Vassilevich {\it Diffeomorphism invariant eigenvalue \break problem
for metric perturbations in a bounded region}, SPbU-IP-95, \break
gr-qc/9509051.} \ref{Milton}{Milton,K.A.: Zero point energy of confined
fermions \prD{22}{1980}{1444}.} \ref{MandS}{Mishchenko,A.V. and
Yu.A.Sitenko: Spectral boundary conditions and index theorem for
two-dimensional manifolds with boundary \aop{218}{1992}{199}.}
\ref{Moss}{Moss,I.G. \cqg{6}{1989}{759}.} \ref{MandP}{Moss,I.G. and
S.J.Poletti: Conformal anomaly on an Einstein space with boundary
\pl{B333}{1994}{326}.} \ref{MandP2}{Moss,I.G. and S.J.Poletti
\np{341}{1990}{155}.} \ref{NandOC}{Nash, C. and O'Connor,D.J.: Determinants
of Laplacians, the Ray-Singer torsion on lens spaces and the Riemann zeta
function \jmp{36}{1995}{1462}.} \ref{NandS}{Niemi,A.J. and G.W.Semenoff:
Index theorem on open infinite manifolds \np {269}{1986}{131}.}
\ref{NandT}{Ninomiya,M. and C.I.Tan: Axial anomaly and index thorem for
manifolds with boundary \np{245}{1985}{199}.} \ref{norlund2}{N\"orlund~N.
E.:M\'emoire sur les polynomes de Bernoulli. \am {4}{1921} {1922}.}
\ref{Poletti}{Poletti,S.J. \pl{B249}{1990}{355}.} \ref{RandT}{Russell,I.H.
and Toms D.J.: Vacuum energy for massive forms in $R^m\times S^N$,
\cqg{4}{1986}{1357}.} \ref{Rayleigh}{Rayleigh,Lord, {\it The Theory of
Sound}, Vols. I and II, 2nd Edn., (MacMillan, London, 1894).}
\ref{RandS}{R\"omer,H. and P.B.Schroer \pl{21}{77}{182}.}
\ref{Schulman}{Schulman.L.S.} \ref{Sneddon}{Sneddon,I.N.,{\it Mixed Boundary
Value Problems in Potential Theory}, (North-Holland. Amsterdam. 1966).}
\ref{Sommerfeld}{Sommerfeld,A. in {\it Die Differential- und
Integralgleichungen der Mechanik und Physik} by Frank,P. and v.Mises,R.,
(Vieweg, Braunschweig, 1935).} \ref{Trautman}{Trautman,A.: Spinors and Dirac
operators on hypersurfaces \jmp{33}{1992}{4011}.} \ref{VdBandG}{van den
Berg,M. and Gilkey,P. \jfa {120}{1994}{48}.}
\ref{Vass}{Vassilevich,D.V.{Vector fields on a disk with mixed boundary
conditions} gr-qc /9404052.} \ref{Voros}{Voros,A.: Spectral functions,
special functions and the Selberg zeta function. \cmp{110}{1987}439.}
\ref{Watson}{Watson,S. 1998 {\it PhD Thesis.}, University of Bristol.}
\ref{Weis}{Weisberger,W.I. \cmp {112}{1987}{633}.} \ref{Ray}{Ray,D.B.:
Reidemeister torsion and the Laplacian on lens spaces \aim{4}{1970}{109}.}
\ref{McandO}{McAvity,D.M. and Osborn,H. Asymptotic expansion of the heat
kernel for generalised boundary conditions \cqg{8}{1991}{1445}.}
\ref{AandE}{Avramidi,I. and Esposito,G. Heat kernel asymptotics with
generalised boundary conditions, hep-th/9701018.} \ref{MandS}{Moss,I.G. and
Silva P.J., Invariant boundary conditions for gauge theories gr-qc/9610023.}
\ref{barv}{Barvinsky,A.O.\pl{195B}{1987}{344}.} \ref{krantz}{Krantz,S.G.
Partial Differential Equations and Complex Analysis (CRC Press, Boca Raton,
1992).} \ref{treves}{Treves,F. Introduction to Pseudodifferential and
Fourier Integral Operators,\break Vol.1, (Plenum Press,New York,1980).}
\ref{treves2}{Treves,F. {\it Basic linear partial differential equations}
(Academic Press, New York, 1975).} \ref{EandS}{Egorov,Yu.V. and Shubin,M.A.
Partial Differential Equations (Springer-Verlag, Berlin,1991).}
\ref{AandS}{Abramowitz,M. and Stegun,I.A. Handbook of Mathematical Functions
(Dover, New York, 1972).} \ref{ACNY}{Abouelsaood,A., Callan,C.G., Nappi,C.R.
and Yost,S.A.\np{280}{1987} {599}.} \ref{BGKE}{Bordag,M., B.Geyer, K.Kirsten
and E.Elizalde, { Zeta function determinant of the Laplace operator on the
D-dimensional ball}, \cmp{179}{1996}{215}.}
  \ref{DandG}{Driscoll,T.A. and Gottlieb, H.P.W., \prE{68}{2003}{016702}.}
  \ref{avramidi2}{Avramidi,I., {\it Math.Phys.Anal.Geom.} {\bf 7}( 2004) 9;
  arXiv : math-ph/0110020.}
  \ref{DGK}{Dowker,J.S., Gilkey,P.B. and Kirsten,K. {\it Int.J.Math.}
  {\bf 12} (2001) 504.}
  \ref{AandK}{Azzam,A. and Kreysig,E. {\it Siam J. Math.Anal.} {\bf13} (1982)
  254.}
  \ref{GandA}{Gustafson,K. and Abe, T {\it Mathematical Intelligencer}
  {\bf20} (1998) 47, 63.}
  \ref{Gustafson}{Gustafson,K. {\it Cont.Math.} {\bf 218} (1998) 432.}
  \ref{Strauss}{Strauss,W.A. {\it Partial Differential Equations}
  (Wiley, New York, 1992).}
  \ref{pockels}{Pockels,F. {\it \"Uber die partielle Differentialgleichung}
  $\De u+k^2u=0$ (Teubner, Leipzig, 1891).}
  \ref{Siegel}{Siegel,W. \np{109}{1976}{244}.}
  \ref{newton}{Newton,I. {\it Phil.Trans.Roy.Soc.} {\bf22} (1701) 838.}
  \ref{bandle}{Bandle,C. {\it Isoperimetric inequalities and applications},
  (Pitman, London, 1980).}
  \ref{CandJ}{Carslaw,H.S. and Jaeger,J.C. {\it Conduction of heat in
  solids} (Clarendon Press, Oxford, 1959).}
  \ref{Grubb2}{Grubb,G. \cmp{215}{2001}{583}.}
  \ref{Grubb3}{Grubb,G. {\it Ann. Global.Anal.Geom.} {\bf 24} (2003) 1.}
  \ref{Grubb4}{Grubb,G. {\it Nucl.Phys.} B (Proc.Suppl.)
  {\bf104} (2002) 71. }
  \ref{Poincare}{Poincar\'e,H. \ajm{12}{1890}{211}.}
  \ref{Frohlich}{Fr\"ohlich,H. \pr {54}{1938}{945}.}
  \ref{MacRobert}{MacRobert,T.M. {\it Spherical Harmonics} (Methuen, London,
  1947).}
  \ref{CVZ}{Cognola,G., Vanzo,L. and Zerbini,S. \jmp{33}{1992}{222}.}
  \ref{GrandS}{Grubb,G. and Seeley,B. {\it Inv.Mat.} {\bf 121} (1995) 481;
  {\it J.Geom.Analysis} {\bf6} (1996) 31.}
  \ref{BandS}{Br\"uning,J. and Seeley,R.\jfa{95}{1991}{255}. }
  \ref{BandMo}{B\"ar,C. and Moroianu, S., {\it Int.J.Math.} {\bf 14} (2003)
  397.}
  \ref{fulling}{Fulling,S.A. \jpa{36}{2003}{6857}.}
  \ref{BoandF}{Bondurant, J.D. and Fulling, S.A. {\it The
  Dirichlet--to--Robin Transform}; arXiv: math--ph/0408054.}
  \ref{Seeley}{Seeley,R. {\it Comm. Partial Diff. Equns.} {\bf 27} (2002)
  2403.}
  \ref{Seeley2}{Seeley,R. {\it Int.J.Mod.Phys.} A {\bf 18} (2003) 2197.}
  \ref{GandG}{Gilkey,P.B. and Grubb,G.}
  \ref{JLNP}{Jakobson,D.,Levitin,M., Nadirashvili,N. and Polterovich, I.,
  {\it Spectral problems with mixed Dirichlet--Neumann boundary conditions:
  isospectrality and beyond} arXiv: math-sp/0409154.}
  \ref{downud}{Dowker,J.S. {\it Nucl.Phys.} B (Proc.Suppl.)
  {\bf104} (2002) 153. }
  \ref{RandS}{Romeo,A. and Saharian,A.A. \jpa{35}{2002}{1297}.}
  \ref{Solod}{Solodukhin,S.N. \prD{63}{2001}{044002}.}
  \ref{AandC}{de Albuquerque,L.C. and Cavalcanti,R.M. \jpa{37}{2004}{7039}.}
  \ref{Kac}{Kac,M. \amm{73}{1966}{1}.}
  \ref{FMPS}{Falomir,H.,Muschietti,M.A., Pisani,P.A.G. and Seeley,R. {\it
  Unusual poles of the \zfs\ for some regular singular differential
  operators}, arXiv: math-ph/0303030.}
  \ref{GandG}{Gilkey,P.B. and Grubb,G. {\it Comm. Partial Diff. Equns.} {\bf
  23} (1998) 777.}
  \ref{CandZ}{Cognola,G. and Zerbini,S. {\it Lett.Math.Phys.} {\bf 42}(1997)
  95.}
  \ref{Kirstenb}{Kirsten,K. {\it Spectral functions in mathematics
  and physics}, (Chapman and Hall/ CRC Press, Boca Raton, 2001). }

\end{putreferences}
\bye